\documentclass[11pt]{article}

\usepackage{amsmath,amsthm,amssymb,amsfonts,graphicx}
\usepackage[all,2cell,ps]{xy}

\theoremstyle{plain}
\newtheorem{thm}{Theorem}[section]
\newtheorem{lem}[thm]{Lemma}
\newtheorem{prop}[thm]{Proposition}
\newtheorem{cor}[thm]{Corollary}

\theoremstyle{definition}

\newtheorem*{rem}{Remark}

\theoremstyle{remark}

\newcommand{\ep}{\epsilon}

\DeclareMathOperator{\U}{U}

\newcommand{\Z}{\mathbb Z}

\newcommand{\R}{\mathbb R}

\newcommand{\Q}{\mathbb Q}

\newcommand{\gam}{\gamma}

\newcommand{\lam}{\lambda}
\newcommand{\Lam}{\Lambda}

\newcommand{\C}{\mathbb C}

\newcommand{\thet}{\theta}
\newcommand{\Gam}{\Gamma}

\newcommand{\conj}{\overline}

\newcommand{\Del}{\Delta}
\newcommand{\mc}{\mathcal}
\newcommand{\del}{\delta}

\DeclareMathOperator{\PSL}{PSL}
\DeclareMathOperator{\GL}{GL}

\DeclareMathOperator{\PU}{PU}
\DeclareMathOperator{\SU}{SU}

\DeclareMathOperator{\real}{Re}

\DeclareMathOperator{\Ad}{Ad}
\DeclareMathOperator{\Tr}{Tr}

\DeclareMathOperator{\Gal}{Gal}

\newenvironment{pf}{\begin{proof}}{\end{proof}}

\newenvironment{enum}{\begin{enumerate}}{\end{enumerate}}

\title{Arithmeticity of complex hyperbolic triangle groups}

\author{Matthew Stover\footnote{partially supported by NSF RTG grant DMS 0602191}\\ University of Michigan\\ \small{\textsf{stoverm@umich.edu}}}

\date{\today}

\begin{document}

\maketitle

\begin{abstract}
Complex hyperbolic triangle groups, originally studied by Mostow in building the first nonarithmetic lattices in $\PU(2, 1)$, are a natural generalization of the classical triangle groups. A theorem of Takeuchi states that there are only finitely many Fuchsian triangle groups that determine an arithmetic lattice in $\PSL_2(\R)$, so triangle groups are generically nonarithmetic. We prove similar finiteness theorems for complex hyperbolic triangle groups that determine an arithmetic lattice in $\PU(2, 1)$.
\end{abstract}


\section{Introduction}\label{intro}


In his seminal 1980 paper, Mostow constructed lattices in $\PU(2, 1)$ generated by three complex reflections \cite{Mostow}. He not only gave a new geometric method for building lattices acting on the complex hyperbolic plane, but gave the first examples of nonarithmetic lattices in $\PU(2, 1)$. Complex reflection groups are a generalization of groups generated by reflections through hyperplanes in constant curvature spaces, and Mostow's groups are a natural extension to the complex hyperbolic plane of the classical triangle groups. They are often called \emph{complex hyperbolic triangle groups}. We introduce these groups in $\S$\ref{triangle groups}. See also \cite{Goldman--Parker, Schwartz}, which, along with \cite{Mostow}, inspired much of the recent surge of activity surrounding these groups.

Around the same time, Takeuchi classified the Fuchsian triangle groups that determine arithmetic lattices in $\PSL_2(\R)$ \cite{Takeuchi}. In particular, he proved that there are finitely many and gave a complete list. Since there are infinitely many triangle groups acting on the hyperbolic plane discretely with finite covolume, triangle groups are generically nonarithmetic. The purpose of this paper is to give analogous finiteness results for complex hyperbolic triangle groups that determine an arithmetic lattice in $\PU(2, 1)$.

A particular difficulty with complex hyperbolic triangle groups is that the complex triangle is not uniquely determined by its angles. One must also consider the so-called \emph{angular invariant} $\psi \in [0, 2 \pi)$. See $\S$\ref{triangle groups}. In particular, there is a $1$-dimensional deformation space of complex triangles with fixed triple of angles. The typical assumption is that $\psi$ is a rational multiple of $\pi$, in which case the angular invariant is called \emph{rational}. We call it \emph{irrational} otherwise.

When a complex hyperbolic triangle group is also an arithmetic lattice, we will call it an arithmetic complex hyperbolic triangle group. Note that this immediately implies discreteness. Our first result is for nonuniform arithmetic complex hyperbolic triangle groups. We prove the following in $\S$\ref{finiteness}.


\begin{thm}\label{intro nonuniform}
There are finitely many nonuniform arithmetic complex hyperbolic triangle groups with rational angular invariant. If $\Gam$ is a nonuniform arithmetic complex hyperbolic triangle group with irrational angular invariant $\psi$, then $e^{i \psi}$ is contained in a biquadratic extension of $\Q$.
\end{thm}


We emphasize that complex reflection groups are allowed to have generators of arbitrary finite order. A usual assumption is that all generators have the same order, a restriction that we avoid. See Theorem \ref{complete nonuniform} for a more precise formulation of Theorem \ref{intro nonuniform}. Proving that a candidate is indeed a lattice is remarkably difficult, as evidenced in \cite{Mostow, Deraux--Parker--Paupert}, so we do not give a definitive list. One consequence of the proof (see Theorem \ref{intro arithmetic restrictions}(1) below) is the following.


\begin{cor}\label{intro nonuniform corollary}
Suppose that $\Gam$ is a nonuniform lattice in $\U(2, 1)$. If $\Gam$ contains a complex reflection of order $5$ or at least $7$, then $\Gam$ is nonarithmetic.
\end{cor}


In the cocompact setting, the arithmetic is much more complicated. Arithmetic subgroups of $\U(2, 1)$ come in two types, defined in $\S$\ref{arithmetic}, often called first and second. In $\S$\ref{proof of first type} we prove the following auxiliary result, generalizing a well-known fact for hyperbolic reflection groups.


\begin{thm}\label{intro fuchsian}
Let $\Gam < \U(2, 1)$ be a lattice containing a complex reflection. Then $\Gam$ contains a Fuchsian subgroup stabilizing the wall of the reflection in $\mathbf H_\C^2$.
\end{thm}


We also give a generalization to higher-dimensional complex reflection groups. Theorem \ref{intro fuchsian} leads to the following, which we also prove in $\S$\ref{proof of first type}.


\begin{thm}\label{intro first type}
Let $\Gam < \U(2, 1)$ be a lattice, and suppose that $\Gam$ is commensurable with a lattice $\Lam$ containing a complex reflection. Then $\Gam$ is either arithmetic of first type or nonarithmetic.
\end{thm}


In particular, when considering a complex reflection group as a candidate for a nonarithmetic lattice, one must only show that it is not of the first type. Fortunately, this is the case where the arithmetic is simplest to understand.

The effect of the angular invariant is a particular sticking point in generalizing Takeuchi's methods. In $\S$\ref{data}, the technical heart of the paper, we study the interdependence between the geometric invariants of the triangle and the arithmetic of the lattice. We collect the most useful of these facts as the following. See $\S\S$\ref{triangle groups}-\ref{arithmetic} for our notation.


\begin{thm}\label{intro arithmetic restrictions}
Suppose that $\Gam$ is an arithmetic complex hyperbolic triangle group. Suppose that for $j = 1, 2, 3$ the generators have reflection factors $\eta_j$, the complex angles of the triangle are $\thet_j$, and that the angular invariant is $\psi$. Let $E$ be the totally imaginary quadratic extension of the totally real field $F$ that defines $\Gam$ as an arithmetic lattice. Then:
\begin{enum}
\item $\eta_j \in E$ for all $j$;
\item $\cos^2 \thet_j \in F$ for all $j$;
\item $e^{2 i \psi} \in E$ and $\cos^2 \psi \in F$;
\item If $\thet_j \leq \pi / 3$ for all $j$, then
\[
\cos^2 \psi \in \Q \big( \cos^2 \thet_1, \cos^2 \thet_2, \cos^2 \thet_3, \cos \thet_1 \cos \thet_2 \cos \thet_3 \big);
\]
\item $E \subseteq \Q \big( \cos^2 \thet_1, \cos^2 \thet_2, \cos^2 \thet_3, e^{i \psi} \cos \thet_1 \cos \thet_2 \cos \thet_3 \big)$;
\item If $\psi$ is rational, then $E$ is a subfield of a cyclotomic field.
\end{enum}
\end{thm}


In $\S$\ref{finiteness}, we use the results from $\S$\ref{data} to prove finiteness results for cocompact arithmetic complex hyperbolic triangle groups with rational angular invariant. We also give restrictions for irrational angular invariants, though it is unknown whether such a lattice exists. When the complex triangle is a right triangle, we prove the following.


\begin{thm}\label{intro right}
Suppose that $\Gam$ is an arithmetic complex hyperbolic triangle group for which the associated complex triangle is a right triangle. Then the angles of the triangle are the angles of an arithmetic Fuchsian triangle group. There are finitely many such $\Gam$ with rational angular invariant.
\end{thm}


Finally, we consider equilateral triangles at the end of $\S$\ref{finiteness}. This is the case which has received the most attention, in particular from Mostow \cite{Mostow} and, in the ideal case, by Goldman--Parker \cite{Goldman--Parker} and Schwartz \cite{Schwartz}. See also \cite{Deraux}. Here we cannot explicitly bound orders of generators, angles, or angular invariants because our proof relies on asymptotic number theory for which we do not know precise constants. Nevertheless, we obtain finiteness in the situation that has received the greatest amount of attention since Mostow's original paper. See \cite{Parker, Parker--Paupert, Paupert, Deraux--Parker--Paupert} and references therein for more recent examples of lattices and restrictions on discreteness.


\begin{thm}\label{intro equilateral}
There are finitely many arithmetic complex hyperbolic equilateral triangle groups with rational angular invariant.
\end{thm}


\section{Complex hyperbolic triangle groups}\label{triangle groups}


We assume some basic knowledge of complex hyperbolic geometry, e.g., the first three chapters of \cite{Goldman}. Let $V$ be a three-dimensional complex vector space, equipped with a hermitian form $h$ of signature $(2, 1)$. Complex hyperbolic space $\mathbf H_\C^2$ is the space of $h$-negative lines in $V$. The metric on $\mathbf H_\C^2$ is defined via $h$ as in \cite[Chapter 3]{Goldman}, and the action of $\U(2, 1)$ on $\mathbf H_\C^2$ by isometries descends from its action on $V$ and factors through projection onto $\PU(2, 1)$. Its ideal boundary $\partial \mathbf H_\C^2$ is the space of $h$-isotropic lines, and we set $\conj{\mathbf H}_\C^2 = \mathbf H_\C^2 \cup \partial \mathbf H_\C^2$.

A \emph{complex reflection} is a diagonalizable linear map $R : V \to V$ with one eigenvalue of multiplicity $2$ (or, more generally, multiplicity $n - 1$ when $\dim(V) = n$). We assume that $R$ has finite order, so the third eigenvalue is a root of unity $\eta$. We call $\eta$ the \emph{reflection factor} of $R$. Decompose $V = V_1 \oplus V_\eta$ into the 1- and $\eta$-eigenspaces, and choose $v_\eta \in V$ so that $V_\eta = \textrm{Span}_\C \{ v_\eta \}$. We begin with an elementary lemma that will be of use later, keeping in mind that every complex reflection has $1$ as an eigenvalue.


\begin{lem}\label{eigenvalues in E}
Let $A \in \GL_n(\C)$ be a diagonalizable linear transformation. Let $E \subseteq \C$ be a subfield, and suppose that $E^n$ has a basis consisting of eigenvectors for $A$. Furthermore, suppose that $A$ has at least one eigenvalue in $E$ and that there exists $x \in \C^\times$ so that $x A \in \GL_n(E)$. Then all eigenvalues of $A$ are in $E$.
\end{lem}


\begin{pf}
Let $v_1, \dots, v_n \in E^n$ be a basis of eigenvectors for $A$, and let $\lam_j$ be the eigenvalue associated with $v_j$, $1 \leq j \leq n$. Without loss of generality, $\lam_1 \in E$. Then $x A$ also has eigenvectors $v_1, \dots, v_n$, and $x A v_j = x \lam_j v_j \in E^n$ for all $j$, since $x A \in \GL_n(E)$. Then $x \lam_j \in E$, $1 \leq j \leq n$. Since $\lam_1 \in E$,  it follows that $x \in E$, which implies that $\lam_j \in E$ for all $j$.
\end{pf}


Assume that $R \in \U(2, 1)$. Then the fixed point set of $R$ acting on $\mathbf H_\C^2$ is the subset of $h$-negative lines in $V_1$. This is a totally geodesic holomorphic embedding of the hyperbolic plane if and only if $V_\eta$ is an $h$-positive line. These subspaces are called \emph{complex hyperbolic lines}. Following \cite[$\S$3.1]{Goldman}, we call $v_\eta$ a \emph{polar vector} for $R$.

When $V_\eta$ is $h$-negative, the fixed set of $R$ on $\mathbf H_\C^2$ is a point, and $R$ is sometimes called a reflection through that point. The complex reflections in this paper will always be of through complex hyperbolic lines. That is, the $\eta$-eigenspace will always be an $h$-positive line.

Let $W$ be the complex hyperbolic line in $\mathbf H_\C^2$ fixed by $R$. We call this the \emph{wall} of $R$. If $v_\eta$ is a polar vector, then $R$ is the linear transformation
\begin{equation}\label{reflection equation}
z \mapsto z + (\eta - 1)\frac{h(z, v_\eta)}{h(v_\eta, v_\eta)} v_\eta.
\end{equation}
We refrain from normalizing the polar vector to have $h$-norm one, since we will often choose a polar vector with coordinates in a subfield $E$ of $\C$, and $E^3 \subset V$ might not contain an $h$-norm one representative for the given line of polar vectors.

Now, consider three complex reflections $R_1, R_2, R_3 \in \U(2, 1)$ with respective distinct walls $W_1, W_2, W_3$ in $\mathbf H_\C^2$. If $v_j$ is a polar vector for $R_j$, then $W_j$ and $W_{j + 1}$ (with cyclic indices) meet in $\mathbf H_\C^2$ if and only if
\begin{equation}\label{cosine equation}
h(W_j, W_{j + 1}) = \frac{|h(v_j, v_{j + 1})|^2}{h(v_j, v_j) h(v_{j + 1}, v_{j + 1})} < 1,
\end{equation}
The two walls meet at a point $z_j$ stabilized by the subgroup of $\U(2, 1)$ generated by $R_j$ and $R_{j + 1}$. The \emph{complex angle} $\thet_j$ between $W_j$ and $W_{j + 1}$, the minimum angle between the two walls, satisfies $\cos^2 \thet_j = h(W_j, W_{j + 1})$.

The walls $W_j$ and $W_{j + 1}$ meet at a point $p_j$ in $\partial \mathbf H_\C^2$ if and only if
\begin{equation}\label{parabolic cosine}
\frac{|h(v_j, v_{j + 1})|^2}{h(v_j, v_j) h(v_{j + 1}, v_{j + 1})} = 1,
\end{equation}
so we say that the complex angle is zero. The group generated by $R_j$ and $R_{j + 1}$ fixes $p_j$, so it is contained in a parabolic subgroup of $\U(2, 1)$. See \cite[$\S$3.3.2]{Goldman}.

Let $\{R_j\}$ be reflections through walls $\{W_j\}$, $j = 1, 2, 3$. When the pairwise intersections of the walls are nontrivial in $\conj{\mathbf H}_\C^2$, they determine a \emph{complex triangle} in $\mathbf H_\C^2$, possibly with ideal vertices. The subgroup $\triangle(R_1, R_2, R_3)$ of $\U(2, 1)$ generated by the $R_j$s is called a \emph{complex hyperbolic triangle group}.

A complex hyperbolic triangle group is sometimes defined as one with order two generators, and groups with higher order generators are called \emph{generalized} triangle groups. We avoid this distinction and do not make the usual assumption that all generators have the same order.

Unlike Fuchsian triangle groups, the complex angles $\{ \thet_1, \thet_2, \thet_3 \}$ do not suffice to determine $\triangle(R_1, R_2, R_3)$ up to $\mathrm{Isom}(\mathbf H_\C^2)$-equivalence. We also need to consider Cartan's \emph{angular invariant}
\begin{equation}\label{angular definition}
\psi = \mathrm{arg} \big( h(v_1, v_2) h(v_2, v_3) h(v_3, v_1) \big).
\end{equation}
A complex triangle is uniquely determined up to complex hyperbolic isometry by the complex angles between the walls, and the angular invariant. See \cite{Brehm} and \cite[Proposition 1]{Pratoussevitch}. Up to the action of complex conjugation on $\mathbf H_\C^2$, we can assume $\psi \in [0, \pi]$.

We call the angular invariant \emph{rational} if $\psi = s \pi /t$ for some (relatively prime) $s, t \in \Z$. In other words, the angular invariant is rational if and only if $e^{i \psi}$ is a root of unity.

Let $\triangle(R_1, R_2, R_3)$ be a complex hyperbolic triangle group in $\U(2, 1)$ with reflection factors $\eta_j$, complex angles $\thet_j$, polar vectors $v_j$, $j = 1, 2, 3$, and angular invariant $\psi$. Suppose that $\{v_1, v_2, v_3\}$ is a basis for $V$. Then $\triangle(R_1, R_2, R_3)$ preserves the hermitian form
\begin{equation}\label{h triangle}
h_{\triangle(R_1, R_2, R_3)} = \begin{pmatrix}
1 &
e^{i \psi} \cos \thet_1 &
e^{i \psi} \cos \thet_3 \\

e^{-i \psi} \cos \thet_1 &
1 &
e^{i \psi} \cos \thet_2 \\

e^{-i \psi} \cos \thet_3 &
e^{-i \psi} \cos \thet_2 &
1
\end{pmatrix}.
\end{equation}
We denote this by $h_\triangle$ when the generators are clear.


\section{Arithmetic subgroups of $\U(2, 1)$}\label{arithmetic}


Let $F$ be a totally real number field, $E$ a totally imaginary quadratic extension, and $\mc D$ a central simple $E$-algebra of degree $d$. Let $\tau : \mc D \to \mc D$ be an involution, that is, an antiautomorphism of order two. Then $\tau$ is of \emph{second kind} if $\tau |_E$ is the Galois involution of $E / F$. There are two cases of interest.
\begin{enum}

\item If $\mc D = E$ (i.e., $d = 1$), then $\tau$ is the Galois involution.

\item If $d = 3$, then $\mc D$ is a cubic division algebra with center $E$.

\end{enum}
See \cite{Knus--et.al.} for more on algebras with involution.

For $d$ as above, let $r = 3 / d$. A form $h : \mc D^r \to \mc D$ is called \emph{hermitian} or $\tau$-\emph{hermitian} if it satisfies the usual definition of a hermitian form with $\tau$ in place of complex conjugation. If $d = 1$, then $h$ is a hermitian form on $E^3$ as usual. If $d = 3$, then there exists an element $x \in \mc D^*$ such that $\tau(x) = x$ and $h(y_1, y_2) = \tau(y_1) x y_2$ for all $y_1, y_2 \in \mc D$.

This determines an algebraic group $\mc G$, the group of elements in $\GL_r(\mc D)$ preserving $h$. For every embedding $\iota : F \to \R$, we obtain an embedding of $\mc G$ into the real Lie group $\U(\iota(h))$. Let $\conj{\mc G}$ be the associated projective unitary group.

If $\mc O$ is a order in $\mc D^r$, then the subgroup $\Gam_{\mc O}$ of $\GL_r(\mc O)$ preserving $h$ embeds as a discrete subgroup of
\[
\mc G(\R) = \prod_{\iota : F \to \R} \U(\iota(h)).
\]
If $\conj{\Gam}_{\mc O}$ is the image of $\Gam_{\mc O}$ in $\conj{\mc G}$, then $\conj{\Gam}_{\mc O}$ is a discrete subgroup of the associated product of projective unitary groups.

The projection of $\Gam_{\mc O}$ onto any factor of $\mc G(\R)$ is discrete if and only if the kernel of the projection of $\mc G(\R)$ onto that factor is compact. Therefore, we obtain a discrete subgroup of $\U(2, 1)$ if and only if $\U(\iota(h))$ is noncompact for exactly one real embedding of $F$.

Then $\conj{\Gam}_{\mc O}$ is a lattice in $\PU(2, 1)$ by the well-known theorem of Borel and Harish-Chandra. An arithmetic lattice in $\PU(2, 1)$ is any lattice $\Gam < \PU(2, 1)$ which is commensurable with $\conj{\Gam}_{\mc O}$ for some $\mc G$ as above and an order $\mc O$ in $\mc D$.

Since arithmeticity only requires commensurability with $\Gam_{\mc O}$, studying an arbitrary $\Gam$ in the commensurability class of $\Gam_{\mc O}$ requires great care. The image of any element $\gam \in \Gam$ in $\PU(2, 1)$ does, however, have a representative in $\GL_3(E)$, that is, if there exists $x \in \C^\times$ so $x \gam \in \GL_3(E)$. This follows from the fact, due to Vinberg \cite{Vinberg}, that $\Gam$ is $F$-defined over the \emph{adjoint form} $\conj{\mc G}$, i.e.,
\[
\Q \big( \Tr \Ad \Gam \big) = F.
\]
This important fact also follows from \cite[Proposition 4.2]{Platonov--Rapinchuk}.


\section{Proofs of Theorems \ref{intro fuchsian} and \ref{intro first type}}\label{proof of first type}


We require some elementary results from the theory of discrete subgroups of Lie groups. The primary reference is \cite{Raghunathan}. Let $G$ be a second countable, locally compact group and $\Gam < G$ a lattice. Recall that $G / \Gam$ carries a finite $G$-invariant measure and $\Gam$ is \emph{uniform} in $G$ if $G / \Gam$ is compact. For a subgroup $H < G$, we let $Z_G(H)$ denote the centralizer of $H$ in $G$. We need the following two results.


\begin{lem}[\cite{Raghunathan} Lemma 1.14]\label{Rag1.14}
Let $G$ be a second countable locally compact group, $\Gam < G$ a lattice, $\Del \subset \Gam$ a finite subset, and $Z_G(\Del)$ the centralizer of $\Del$ in $G$. Then, $Z_G(\Del) \Gam$ is closed in $G$.
\end{lem}


\begin{thm}[\cite{Raghunathan} Theorem 1.13]\label{Rag1.13}
Let $G$ be a second countable locally compact group, $\Gam < G$ be a uniform lattice, and $H < G$ be a closed subgroup. Then $H \Gam$ is closed in $G$ if and only if $H \cap \Gam$ is a lattice in $H$.
\end{thm}


We now use the above to prove Theorem \ref{intro fuchsian}.


\begin{pf}[Proof of Theorem \ref{intro fuchsian}]
Assume that $\Gam$ is a cocompact arithmetic lattice in $\U(2, 1)$ containing a complex reflection and that $\Del$ is the subgroup of $\Gam$ generated by this reflection. The centralizer of $\Del$ in $\U(2, 1)$ is isomorphic to the extension of $\U(1, 1)$ by the center of $\U(2, 1)$, and is the stabilizer in $\U(2, 1)$ of the wall of the reflection that generates $\Del$. It follows from Lemma \ref{Rag1.14} and Theorem \ref{Rag1.13} that $\Gam \cap \U(1, 1)$ is a lattice. Since any sublattice of an arithmetic lattice is arithmetic, $\Gam$ contains a totally geodesic arithmetic Fuchsian subgroup.
\end{pf}


\begin{pf}[Proof of Theorem \ref{intro first type}]
A totally geodesic arithmetic Fuchsian group comes from a subalgebra of $\mc D^r$, with notation as in $\S$\ref{arithmetic}. When $\Gam$ is of second type, $\mc D$ is a cubic division algebra. The totally geodesic Fuchsian group would correspond to a quaternion subalgebra of $\mc D$, which is impossible. When $\Gam$ is of first type, this quaternion subalgebra corresponds to rank 2 subspaces of $E^3$ on which $h$ has signature (1, 1). Therefore, $\Gam$ contains complex reflections if and only if $\Gam$ is of first type.
\end{pf}


\begin{rem}
One can also prove Theorem \ref{intro first type} using the structure of unit groups of division algebras.
\end{rem}


We now briefly describe how these results generalize to reflections acting on higher-dimensional complex hyperbolic spaces. If $\Gam < \U(n, 1)$ is a lattice, an element $R \in \Gam$ is a \emph{codimension} $s$ \emph{reflection} if it stabilizes a totally geodesic embedded $\mathbf H_\C^{n - s}$ and acts by an element of the unitary group of the normal bundle to the wall. If $\Gam$ is arithmetic, the associated algebraic group is constructed via a hermitian form on $\mc D^r$, where $\mc D$ is a division algebra of degree $d$ with involution of the second kind over a totally imaginary field $E$, and where $r d = n + 1$.


\begin{thm}\label{generalized reflections}
Suppose $\Gam < \U(n, 1)$ is a cocompact arithmetic lattice with associated algebraic group coming from a hermitian form on $\mc D^r$, where $\mc D$ is a central simple algebra with involution of the second kind. If $\Gam$ contains a codimension $s$ reflection, then $\Gam$ contains a cocompact lattice in $\U(n - s, 1)$. Also, $n - s + 1 = \ell d$ for some $1 < \ell \leq r$ and the associated algebraic subgroup comes from a hermitian form on $\mc D^\ell$.
\end{thm}


\begin{cor}\label{higher dimensional reflection groups}
Let $\Gam < \U(n, 1)$ be an arithmetic lattice generated by complex reflections through totally geodesic walls isometric to $\mathbf H_\C^{n - 1}$. Then $\Gam$ is of so-called first type, i.e., the associated algebraic group is the automorphism group of a hermitian form on $E^{n + 1}$, where $E$ is some totally imaginary quadratic extension of a totally real field.
\end{cor}


\section{Arithmetic data for complex hyperbolic triangle groups}\label{data}


In this section, we relate the geometric invariants of a complex triangle to the arithmetic invariants of the complex reflection group. It is the technical heart of the paper.

Let $\Gam = \triangle(R_1, R_2, R_3)$ be a complex hyperbolic triangle group with reflection factors $\eta_j$, complex angles $\thet_j$, and angular invariant $\psi$. Assume that $\Gam$ is an arithmetic lattice in $\U(2, 1)$. By Theorem \ref{intro first type}, $\Gam$ is of first type, so there is an associated hermitian form $h$ over a totally imaginary field $E$. Let $F$ be the totally real quadratic subfield of $E$.


\begin{lem}\label{polar vectors in E}
We can choose polar vectors $v_j$ for the reflection $R_j$ so that $v_j \in E^3$.
\end{lem}


\begin{pf}
Associated with each reflection is an arithmetic Fuchsian subgroup of $\Gam$. When $\Gam$ is a uniform lattice, this follows from Theorem \ref{intro fuchsian}. For the nonuniform case, see \cite[Chapter 5]{Holzapfel}. Arithmetic Fuchsian subgroups stabilizing a complex hyperbolic line come from the $h$-orthogonal complement of an $h$-positive line in $E^3$. (To be more precise, this line is $h$-positive over the unique real embedding of $F$ at which $h$ is indefinite.) Any vector in $E^3$ representing this line is a polar vector for $R_j$.
\end{pf}


This leads us to the following important fact.


\begin{lem}\label{reflection factors in E}
Each reflection factor $\eta_j$ is contained in $E$.
\end{lem}


\begin{pf}
It follows from Proposition 4.2 in \cite{Platonov--Rapinchuk} that there exists an $x_j \in \C$ so that $x_j R_j \in \GL_3(E)$ (see the end of $\S$\ref{arithmetic} above). By Lemma \ref{polar vectors in E}, and because the $h$-orthogonal complement to a polar vector evidently has an $E$-basis, $E^3$ has a basis of eigenvectors for $R_j$. The lemma follows from Lemma \ref{eigenvalues in E}.
\end{pf}


Now we turn to the complex angles and the angular invariant.


\begin{lem}\label{angle squared in E}
For each $j$, $\cos^2 \thet_j \in F$ and $e^{2 i \psi} \in E$.
\end{lem}


\begin{pf}
Choose polar vectors $v_j \in E^3$. The terms in Equations \eqref{cosine equation} and \eqref{parabolic cosine} resulting from these choices of polar vectors are all contained in $E$. Hence $\cos^2 \thet_j \in F$. One can also prove this using $\Tr \Ad(R_1 R_2)$ and Lemma \ref{reflection factors in E}.

Similarly, consider
\[
\del = h(v_1, v_2) h(v_2, v_3) h(v_3, v_1) = r e^{i \psi} \in E
\]
from Equation \eqref{angular definition}. Note that $e^{i \psi} \in E$ if and only if $r \in E$. Either way,
when $\del \neq 0$, we have $\del / \conj \del = e^{2 i \psi} \in E$. This completes the proof.
\end{pf}


Combining the above, we see that
\[
\Q \big( \eta_1, \eta_2, \eta_3, \cos^2 \thet_1, \cos^2 \thet_2, \cos^2 \thet_3, e^{2 i \psi} \big) \subseteq E.
\]
We can also bound $E$ from above using the fact that $E \subseteq \Q \big( \Tr \Gam \big)$. Using well-known computations of traces for products of reflections (e.g., \cite[$\S$4]{Mostow} or \cite{Pratoussevitch}), we have
\[
\Q \big( \Tr \Gam \big) = \Q(\eta_1, \eta_2, \eta_3, \cos^2 \thet_1, \cos^2 \thet_2, \cos^2 \thet_3, e^{i \psi} \cos \thet_1 \cos \thet_2 \cos \thet_3 \big).
\]
Similarly,
\[
\Q \big( \real(\eta_1), \real(\eta_2), \real(\eta_3), \cos^2 \thet_1, \cos^2 \thet_2, \cos^2 \thet_3, \cos^2 \psi \big) \subseteq F \subseteq
\]
\[
\Q(\real(\eta_1), \real(\eta_2), \real(\eta_3), \cos^2 \thet_1, \cos^2 \thet_2, \cos^2 \thet_3, \cos \psi \cos \thet_1 \cos \thet_2 \cos \thet_3 \big).
\]
This gives the following.


\begin{cor}\label{E cyclotomic}
Let $\Gam$ be a complex hyperbolic triangle group and an arithmetic lattice in $\U(2, 1)$. If the angular invariant of the triangle associated with $\Gam$ is rational, then the fields that define $\Gam$ as an arithmetic lattice are subfields of a cyclotomic field.
\end{cor}


Let $h_\triangle$ be as in \eqref{h triangle} and consider $h_\triangle$ as a hermitian form on the extension
\[
E_\triangle = \Q \big( \eta_1, \eta_2, \eta_3, \cos \thet_1, \cos \thet_2, \cos \thet_3, e^{i \psi} \big),
\]
of $E$. It follows from \cite[$\S$2]{Mostow} that $h$ and $h_\triangle$ are equivalent over $E_\triangle$. Consequently, $h_\triangle$ is indefinite over exactly one complex conjugate pair of places of $E$. This implies that there are precisely $[E_\triangle : E]$ conjugate pairs of places of $E_\triangle$ over which $h_\triangle$ is indefinite.

Let $H$ be a hermitian in $3$ variables over the complex numbers for which there is a vector with positive $H$-norm. Then $H$ is indefinite if and only if $\det(H) < 0$. Since any polar vector has positive $h_\triangle$-norm by definition, we have the following.


\begin{prop}\label{indefinite places}
There are exactly $[E_\triangle : E]$ complex conjugate pairs of Galois automorphisms $\tau$ of $E_\triangle \subset \C$ under which $\tau \left( \det(h_\triangle) \right)$ is negative. All such automorphisms act trivially on $E$.
\end{prop}


This has the following consequence for the relationship between the geometry of the triangle and the arithmetic of the lattice.


\begin{cor}\label{reflection factors are small}
If $\Gam$ is a complex hyperbolic triangle group and an arithmetic lattice, then the reflection factors of $\Gam$ are restricted by the geometry of the triangle. In particular,
\[
E_\triangle = \Q \big( \cos \thet_1, \cos \thet_2, \cos \thet_3, e^{i \psi} \big).
\]
\end{cor}


\begin{pf}
Since $\det(h_\triangle)$ is independent of the reflection factors, for each Galois automorphism of
\[
E_\triangle / \Q \big( \cos \thet_1, \cos \thet_2, \cos \thet_3, e^{i \psi} \big)
\]
we obtain a new complex conjugate pair of embeddings of $E_\triangle$ into $\C$ such that $\det(h_\triangle)$ is negative. Any such automorphism necessarily acts nontrivially on some reflection factor $\eta_j$. These embeddings of $E_\triangle$ lie over different places of $E$ by Lemma \ref{reflection factors in E}. This contradicts Proposition \ref{indefinite places}.
\end{pf}


We also obtain the following dependence between the angular invariant and the angles of the triangle.


\begin{prop}\label{angular invariant constraint}
If $\Gam$ is a complex hyperbolic triangle group and an arithmetic lattice. If $\Gam$ has rational angular invariant and $\thet_j \leq \pi / 3$ for $j = 1, 2, 3$, then
\[
\cos^2 \psi \in F' = \Q \big( \cos^2 \thet_1, \cos^2 \thet_2, \cos^2 \thet_3, \cos \thet_1 \cos \thet_2 \cos \thet_3 \big).
\]
\end{prop}


\begin{pf}
If $\psi$ is rational, then $E_\triangle$ is a subfield of a cyclotomic field $K_N = \Q \big( \zeta_N \big)$, where $\zeta_N$ is a primitive $N^{th}$ root of unity. Therefore the Galois automorphisms of $E_\triangle$ are induced by $\zeta_N \mapsto \zeta_N^m$ for some $m$ relatively prime to $N$.

Consider the stabilizer $S$ of $F'$ in $\Gal(K_N / \Q)$. It acts on the roots of unity in $E_\triangle$ as a group of rotations along with complex conjugation. By definition of $E_\triangle$, every nontrivial element of $S$ acts nontrivially on $e^{i \psi}$. In particular, if $\cos^2 \psi \notin \Q$ and $S$ contains a rotation through an angle other than an integral multiple of $\pi$, then the orbit of $e^{i \psi}$ under $S$ contains two non-complex conjugate points with distinct negative real parts.

Let $\tau$ be any such automorphism of $E_\triangle$. Then, since $\tau(\cos \thet_j) = \cos \thet_j$ for all $j$ by definition of $S$,
\[
\tau(\det(h_\triangle)) = 1 - \sum_{j = 1}^3 \cos^2 \thet_j + 2 \tau(\cos \psi) \prod_{j = 1}^3 \cos \thet_j.
\]
Furthermore, $1 - \sum \cos^2 \thet_j \leq 0$ for any triple of angles $\thet_j = \pi / r_j$ that are the angles of a hyperbolic triangle with each $r_j \geq 3$. Since $\tau(\cos \psi) < 0$ and $\cos \thet_j > 0$, it follows that $\tau(\det h_\triangle) < 0$. Since $\tau$ acts nontrivially on $e^{2 i \psi} \in E$, this contradicts Proposition \ref{indefinite places}. Therefore, $S$ is generated by complex conjugation and rotation by $\pi$, so $\cos^2 \psi \in F'$.
\end{pf}


\begin{rem}
For several of Mostow's lattices in \cite{Mostow}, $F' = F$ (with notation as above) and $\cos \psi \notin F'$. Thus Proposition \ref{angular invariant constraint} is the strongest possible constraint on rational angular invariants.
\end{rem}


\section{Finiteness results}\label{finiteness}


We are now prepared to collect facts from $\S$\ref{data} to prove Theorem \ref{intro nonuniform}. A more precise version is the following.


\begin{thm}\label{complete nonuniform}
Suppose that $\Gam$ is a complex hyperbolic triangle group and a nonuniform arithmetic lattice in $\U(2, 1)$. Then
\begin{enum}
\item Each generator has order $2$, $3$, $4$, or $6$.
\item Each complex angle $\thet_j$ of the triangle comes from the set
\[
\{\pi / 2, \pi / 3, \pi / 4, \pi / 6, 0\}.
\]
\item If $\psi$ is the angular invariant, then $e^{i \psi}$ lies in a biquadratic extension of $\Q$.
\item If $\psi$ is rational, then $\psi = s \pi / t$ for
\[
t \in \{2, 3, 4, 6, 8, 12\}.
\]
\end{enum}
\end{thm}


\begin{pf}
Since $\Gam$ is a nonuniform arithmetic lattice, the associated field $E$ is imaginary quadratic. For \emph{1}, we apply Lemma \ref{reflection factors in E} to $E$. For \emph{2} and \emph{3}, we apply Lemma \ref{angle squared in E}. Then \emph{4} follows from determining those integers $m$ so that $\varphi(m) = 2$ or $4$ and $e^{2 i \psi}$ is at most quadratic over $\Q$, where $\varphi$ is Euler's totient function.
\end{pf}


See \cite{Paupert, Deraux--Parker--Paupert} for the known nonuniform arithmetic complex hyperbolic triangle groups. We now determine the right triangle groups that can determine an arithmetic lattice in $\SU(2, 1)$.

\begin{pf}[Proof of Theorem \ref{intro right}]
Suppose that $\Gam$ is an arithmetic complex hyperbolic triangle group with $\thet_1 = \pi / 2$. The hermitian form $h_\triangle$ associated with the triangle has determinant
\[
1 - \cos^2 \thet_2 - \cos^2 \thet_3.
\]
By Lemma \ref{angle squared in E}, this is an element of the totally real field $F$ that defines $\Gam$ as an arithmetic lattice. Consequently, there is no Galois automorphism of $F$ over $\Q$ under which this expression remains negative.

This is precisely Takeuchi's condition that determines whether or not the triangle in the hyperbolic plane with angles $\pi / 2, \thet_2, \thet_3$ determines an arithmetic Fuchsian group. The theorem follows from Takeuchi's classification of arithmetic Fuchsian right triangle groups, Lemma \ref{angle squared in E}, and Corollary \ref{reflection factors are small}
\end{pf}


There are $41$ such right triangles in $\mathbf H^2$. We now finish the paper with finiteness for arithmetic complex hyperbolic triangle groups with equilateral complex triangle and rational angular invariant.


\begin{pf}[Proof of Theorem \ref{intro equilateral}]
Let $\Gam$ be an arithmetic complex hyperbolic triangle group with equilateral triangle of angles $\pi / n$ and angular invariant $\psi$. By Proposition \ref{angular invariant constraint}, we can assume that $\psi = s \pi / 12 n$ for some integer $s$. Indeed, $F' = \Q(\cos \pi / n)$, and the assertion follows from an easy Galois theory computation.

Then
\begin{equation}\label{equilateral det}
\det(h_\triangle) = 1 - 3 \cos^2(\pi / n) + 2 \cos(s \pi / 12 n) \cos^3(\pi / n),
\end{equation}
so we want to find a nontrivial Galois automorphism of $F_\triangle$ whose restriction to $F$ is nontrivial and such that the image of \eqref{equilateral det} under this automorphism is negative. Let $p$ be the smallest rational prime not dividing $12 n$. This determines a nontrivial Galois automorphism $\tau_p$ of $F_\triangle$ under which
\begin{equation}\label{equilateral automorphism}
\tau_p(\det(h_\triangle)) = 1 - 3 \cos^2(p \pi / n) + 2 \cos(p s \pi / 12 n) \cos^3(p \pi / n).
\end{equation}
It is nontrivial on $F$ by definition. If we show that $\tau_p(\det(h_\triangle)) < 0$ for $n$ sufficiently large, this, along with Corollary \ref{reflection factors are small}, suffices to prove the theorem.

First, notice that the function
\[
D(x, y) = 1 - 3 \cos^3 + 2 \cos y \cos^3 x
\]
is an increasing function of $x \in (0, \pi / 2)$ for any fixed $y$. In our language, this implies that if $y$ is the angular invariant of an equilateral complex triangle in $\mathbf H_\C^2$ with angle $x$, then it remains an angular invariant for a complex triangle with angle $x'$ for any $x' < x$. Similarly, if we know that $\pi / 12 n$ is an angular invariant for a triangle with angles $p \pi / n$, then we know that $p s \pi / n$ (more precisely, a representative modulo $2 \pi$) is the angular invariant of an equilateral triangle in $\mathbf H_\C^2$ with angles $p \pi / n$. Therefore, it is enough to show that $\pi / 12 n$ is the angular invariant of a triangle having angles $p \pi / n$ for all sufficiently large $n$, where $p$ is the smallest not prime dividing $12 n$.

From the above, we conclude further that it suffices to show that there exists a function $q(n)$ such that $p < q(n)$ and
\begin{equation}\label{equilateral determinant bound}
1 - 3 \cos^2(q(n) \pi / n) + 2 \cos(\pi / 12 n) \cos^3(q(n) \pi / n) < 0
\end{equation}
for all sufficiently large $n$. To prove this, we consider the function $j(n)$, defined by Jacobsthal \cite{Jacobsthal}. For any integer $n$, $j(n)$ is the smallest integer such that any $j(n)$ consecutive integers must contain one that is relatively prime to $n$. Clearly $p \leq j(12 n)$.

Iwaniec \cite{Iwaniec} proved that
\[
j(n) \ll (\log n)^2.
\]
Therefore, for any $\ep > 0$, there is an $n_\ep$ so that the first prime number coprime to $12 n$ is at most $\log(12 n)^{2 + \ep}$ for every $n \geq n_\ep$. Consider the function
\[
f_\ep(x) = 1 - 3 \cos^2 \big( \log(12 / x)^{2 + \ep} \pi x \big) + 2 \cos \big( \pi x / 12 \big) \cos^3 \big( \log(12 / x)^{2 + \ep} \pi x \big).
\]
Then $\lim_{x \to 0} f_\ep(x)$ exists and equals $0$ for all $\ep > 0$. Further, $x = 0$ is a local maximum of $f_\ep$, so $f_\ep(1 / n) < 0$ for all sufficiently large $n$.

Taking $q(n) = \log(n)^{2 + \ep}$ for any small $\ep$ shows that \eqref{equilateral determinant bound} holds for all sufficiently large $n$. This implies that \eqref{equilateral automorphism} is negative for all large $n$. This proves the theorem.
\end{pf}


Unfortunately, the proof of Theorem \ref{intro equilateral} isn't effective, so we cannot list the angles that can possibly determine an arithmetic lattice. In particular, we don't know which $n$ makes the bound from \cite{Iwaniec} effective for any $\ep > 0$. If this bound is less than $n = 10^5$ for some $\ep$, which computer experiments show is extraordinarily likely, then we obtain $n < 5,000,000$. We expect the actual bound to be quite a bit smaller, especially given that the smallest equilateral triangle in $\mathbf H^2$ that defines an arithmetic Fuchsian group has angles $\pi / 15$.


\subsection*{Acknowlegments}

I thank the referee for several helpful comments.


\bibliography{Arithmetic_triangles}


\end{document}